\documentclass[reqno,12pt]{amsart}

\usepackage{mathrsfs}
\usepackage{scalerel}

\theoremstyle{definition}

\theoremstyle{remark}

\usepackage[usenames]{color}

\usepackage{cite}

\usepackage{tikz}

\usepackage[cp1250]{inputenc} % nastavuje pouzite kodovani
\usepackage{a4wide} % nastavuje standardní evropský formát stránek A4
\usepackage[stable]{footmisc} % umoznuje pridavani footnotes i do nazvu kapitol a sekci
\usepackage{fancyhdr} % zahlavi stranek, zapina se \pagestyle{fancy}
\usepackage{layout} % umoznuje prikaz \layout pro zobrazeni diagramu layoutu (nakresleni marginu, paddingu, atd...)

\usepackage{amsthm,amsmath,amstext,amsbsy,amsfonts,amssymb} % AMS
\usepackage{stmaryrd} % extra matematicke symboly, doplnek k AMS
\usepackage{turnstile} % snadne vlastni kresleni symbolu typu dokazovatka, pravdivitka, forsitka, a mnoha dalsich...
\usepackage{mathtools} % rozsireni amsmath
\usepackage{graphicx} % nezbytné pro standardní vkládání obrázků do dokumentu
\usepackage[all,cmtip]{xy} % kresleni diagramu xymatrix, Xy-pic
\usepackage{fancybox} % umožňuje pokročilé rámečkování :-)
\usepackage{hhline} % kresleni vertikalnich i svislych tabulkovych car/dvojcar
\usepackage{colortbl} % barevne tabulky
\usepackage{color} % barevny text
\usepackage{pstricks} % PStricks - kresleni obrazku
\usepackage{qtree} % jednoduche kresleni stromu
\usepackage{tikz} % kresleni komplikovanych diagramu
\usepackage{enumerate} % rozsirene moznosti enumeraci
\usepackage{hyperref} % hypertextove odkazy
\usepackage{index} % nutno pouľít v případě tvorby rejstříku balíčkem makeindex
\usepackage{longtable} % vicestrankove tabulky
\usepackage{cite} % bibliografie

\usepackage{ifthen} % if, then, else
\usepackage{xstring} % prace s retezci znaku
\usepackage{calc} % zakladni aritmetika intuitivne
\usepackage{forloop} % implementace for cyklu

\usepackage{xspace} %\xspace se dava na konec makra, ktere bude pouzivane prevazne v textu, prida pak za jeho uziti mezeru
\ifx\uplodadedmylogic\undefined
	\def\uplodadedmylogic{}
\newcommand{\str}[1]{\StrLen{#1}[\tmplength]\ifthenelse{1=\tmplength{}}{\mathcal{#1}}{\underline{#1}}}

\newcommand{\stru}[2]{\lara{#1,#2}} % struktura
\newcommand{\thn}[1]{\mathrm{#1}} %theory name

\newcommand{\DeLOs}[2]{{\tiny\ifthenelse{\equal{#1}{#2}}{\thn{DeLO}^{#1}}{\ifthenelse{\equal{#2}{}}{\thn{DeLO}^{#1}}{\thn{DeL0}^{\vpair{$#1$}{$#2$}}}}}}

\newcommand{\PAK}{\Rightarrow}
\newcommand{\KAP}{\Leftarrow}

\newcommand{\landef}[2][{}]{\csname Ln#2\endcsname{#1} = \csname L#2\endcsname{#1}}
\newcommand{\landefl}[2][{}]{$\landef[#1]{#2}$, kde \csname Ld#2\endcsname{#1}}

\def\LnfGrpA/{aditivní jazyk teorie grup}
\def\LnFGrpA/{Aditivní jazyk teorie grup}
\def\LnfOrd/{jazyk teorie uspořádání}
\def\LnFOrd/{Jazyk teorie uspořádání}
\def\LnfPrv/{výrokový jazyk s množinou prvovýroků $\Prv$}
\def\LnFPrv/{Výrokový jazyk s množinou prvovýroků $\Prv$}
\def\LnfAr/{jazyk aritmetiky}
\def\LnFAr/{Jazyk aritmetiky}
\def\LnfArZ/{jazyk $\Z$-aritmetiky}
\def\LnFArZ/{Jazyk $\Z$-aritmetiky}
\def\LnfArPA/{jazyk Peanovy aritmetiky}
\def\LnFArPA/{Jazyk Peanovy aritmetiky}
\def\LnfArZPA/{jazyk $\Z$-Peanovy aritmetiky}
\def\LnFArZPA/{Jazyk $\Z$-Peanovy aritmetiky}
\def\LnfArA/{aditivní jazyk aritmetiky}
\def\LnFArA/{Aditivní jazyk aritmetiky}
\def\LnfArAA/{jazyk aditivní aritmetiky}
\def\LnFArAA/{Jazyk aditivní aritmetiky}
\def\LnfArZAA/{jazyk $\Z$-aditivní aritmetiky}
\def\LnFArZAA/{Jazyk $\Z$-aditivní aritmetiky}
\def\LnfArLA/{jazyk lineární aritmetiky}
\def\LnFArLA/{Jazyk lineární aritmetiky}
\def\LnfArZLA/{jazyk $\Z$-lineární aritmetiky}
\def\LnFArZLA/{Jazyk $\Z$-lineární aritmetiky}
\newcommand{\scal}[1]{\underline{#1}} % skalár
\def\LnfArKLA/{jazyk $\kappa$-lineární aritmetiky}
\def\LnFArKLA/{Jazyk $\kappa$-lineární aritmetiky}
\def\LnfArKZLA/{jazyk $\kappa$-$\Z$-lineární aritmetiky}
\def\LnFArKZLA/{Jazyk $\kappa$-$\Z$-lineární aritmetiky}
\def\LnfMod/{jazyk modulů nad okruhem}
\def\LnFMod/{Jazyk modulů nad okruhem}
\def\LnfSuc/{jazyk následníka}
\def\LnFSuc/{Jazyk následníka}
\def\LnfSucO/{jazyk následníka s nulou}
\def\LnFSucO/{Jazyk následníka s nulou}
\newcommand{\Prv}{\mathbb{P}}
\def\DNF/{DNF}
\def\CNF/{CNF}

\fi

\ifx\uplodadedmygeneral\undefined
	\def\uplodadedmygeneral{}
\newcommand{\N}{\mathbb{N}}
\newcommand{\Z}{\mathbb{Z}}

\newcommand{\R}{\mathbb{R}}

\newcommand{\lr}[1]{\{#1\}}

\newcommand{\lara}[1]{\langle #1\rangle}

\newcommand{\set}[2]{\lr{#1;#2}}

\newcommand{\sdiff}{-}
\newcommand{\suni}{\cup}

\newcommand{\vpairg}[3]{{\newbox\horni \newbox\dolni \setbox\horni=\hbox{#1} \setbox\dolni = \hbox{#2} \newdimen\zdeposun \newdimen\wddolni \newdimen\wdhorni \setlength{\wddolni}{\wd\dolni} \setlength{\wdhorni}{\wd\horni} \setlength{\zdeposun}{-\wdhorni-\wddolni/2+\wdhorni/2} {\lower -#3 \hbox{#1}}\kern\zdeposun {\lower #3 \hbox{#2}}}} 
\newcommand{\vpair}[2]{\vpairg{#1}{#2}{.85ex}}

\fi
\ifx\uplodadedmytextstructuring\undefined
	\def\uplodadedmytextstructuring{}
\newcommand{\idf}[3]{\label{#1}\csname index#2\endcsname[|textbf]#3} % in-line definition \idf{label}{name}{definition}
\newcommand{\benum}{\begin{enumerate}}
\newcommand{\eenum}{\end{enumerate}}

\newcommand{\bitem}{\begin{itemize}}
\newcommand{\eitem}{\end{itemize}}

\newcommand{\uvz}[1]{``#1''} % anglicke uvozovky
\newcommand{\axiominfo}[7]{
\begin{tabular}{#7}
Znění: & #1\\
Význam: & #2
\ifthenelse{\equal{#4}{}}{}{\\#3 & #4}%
\ifthenelse{\equal{#6}{}}{}{\\#5 & #6}%
\end{tabular}
}

\newcommand{\repeatstatement}[4][thm]% #1= typ statementu (thm,prop,...) #2 = jmeno labelu, #3=nazev tvrzeni, #4=zneni tvrzeni
{
{
\renewcommand{\thestatement}{\ref{#2}}
\begin{#1}[#3]
#4
\end{#1}
\addtocounter{statement}{-1}
}
}

\newcounter{thislevel}
\newcommand{\mysect}[2]{\setcounter{thislevel}{\value{mysectlevel}} \addtocounter{thislevel}{#1}
\if\arabic{thislevel}0\part{#2}\fi
\if\arabic{thislevel}1\chapter{#2}\fi
\if\arabic{thislevel}2\section{#2}\fi
\if\arabic{thislevel}3\subsection{#2}\fi
\if\arabic{thislevel}4\subsubsection{#2}\fi
\if\arabic{thislevel}5\paragraph{#2}\fi
\if\arabic{thislevel}6\subparagraph{#2}\fi}
\fi
\ifx\uplodadedmypeanoproductsandLA\undefined
	\def\uplodadedmypeanoproductsandLA{}
\newcommand{\dvf}[1]{{#1}^{-1}}

\def\LnFDT/{Language of D-$\Rcal$-theory}

\newcommand{\LdDT}[1]{$0,1$ are constant symbols, $+$ binary functional symbol, $-$ unary functional symbol, $\leq$ binary predicate symbol, $\scal{r}$ unary functional symbols}

\def\LnFDTg/{Language of graded D-$\Rcal$-theory}

\newcommand{\LdDTg}[1]{$0,1$ are constant symbols, $+$ binary functional symbol, $-$ unary functional symbol, $\leq$ binary predicate symbol, $\scal{r}, \dvf{q}$ unary functional symbols}

\newcommand{\Rcal}{\mathcal{R}}

\fi

\begin{document}

%\title{Some remarks on possible values of Peano multiplications}
\title{Metaphors in teaching infinity: Limits, cardinalities, and nonstandard models} %at high school and undergraduate levels
\author{\textsc{Petr~Glivick\' y}}
\address{\textsc{{Petr~Glivick\' y}}: University of Economics, Department of Mathematics\\  %\textsc{\hphantom{Petr~Glivick\' y}}:
Ekonomick\' a 957, 148 00 Praha~4, Czech Republic}
\email{petrglivicky@gmail.com}

\thanks{This paper was processed with contribution of long term institutional support of research activities by Faculty of Informatics and Statistics, University of Economics, Prague.}

%\author{\textsc{Josef~Ml\v cek}}
%\address{\textsc{{Josef~Ml\v cek}}: Charles University, Faculty of Mathematics and Physics, Department of Theoretical Computer Science and Mathematical Logic \\ 
%Malostransk\'e n\'am\v est\' i\ 12, 118 00 Praha~1, Czech Republic}
%\email{josef.mlcek@mff.cuni.cz}

%%    General info
%\subjclass[2010]{03F30, 11U10, 03H15, 03C62}

%\date{March 1, 2013}

%\keywords{Lattice-ordered abelian group, semifield, finitely generated, order-unit}

%\begin{abstract}

%\end{abstract}

% dukazove obraty
\newcommand{\ZNA}[2]{#1$\PAK$#2} % oznaceni kroku dukazu implikace, napr. 1)=>2), ci jen =>
\newcommand{\NAZ}[2]{#1$\KAP$#2} % oznaceni kroku dukazu obracene implikace, napr. 1)<=2), ci jen , <=
\newcommand{\starr}{{}^{*}}

%%%%%%%%%%%%%%%%%%%

\begin{abstract}
Mathematical conception of infinite quantities forms a cornerstone of many disciplines of modern mathematics --- from differential calculus to set theory. In fact, it could be argued that the most significant revolutions in mathematics in the modern period were always triggered by a development in our understanding of infinity. From the pedagogical point of view, the students' comprehension of the concept of infinity is a competence of interdisciplinary value, helping them to grasp the reasons why and how individual disciplines of modern mathematics were built up.

In this paper we present a number of illustrations, examples, and allegories that illuminate and clarify different aspects of infinity. The author has been using and gradually developing these metaphors for infinity in his undergraduate and graduate courses and popular lectures on mathematical logic, set theory, calculus, and nonstandard analysis during the last decade. The aim of this paper is to share the accumulated didactic know-how.
\end{abstract}

\maketitle

\section{Introduction}

The concept of infinity appears in various forms in many disciplines of higher mathematics. For the students new to the mathematical conception of \uvz{infinite quantity}, to internalize the concept and to thoroughly understand all the aspects, nuances, and paradoxes connected with infinity, is often a crucial competence allowing them to gain a deeper insight into subjects as distinct as differential calculus, set theory or Ramsey combinatorics. 

Moreover, in the same way as understanding the abstract concept of geometrical point allows the mind to start discovering the ideal world of classical geometry, comprehending infinity opens the gates to the realms of cardinal and ordinal arithmetics, infinitesimals or existential proofs that are undoubtedly of no less significance for the modern mathematical description of reality than ideal lines and circles were for the ancients.

In what follows, we give a list of examples that illustrate or motivate different aspects of infinity in modern mathematics. Some of the examples are rather classical (just with minimal adjustments), other are completely original. We classify these examples into sections by the aspect of infinity they serve to clarify. We cover three such aspects --- infinity as a limit (Section \ref{sect:limit}), infinite cardinalities and ordinalities (Section \ref{sect:card}), and infinity in nonstandard models of arithmetics and in nonstandard analysis (Section \ref{sect:nonstandard}).

The reader interested in how students perceive these aspects of infinity intuitively prior to their exposition to formal mathematical theories, as well as in the philosophical and historical development of thinking about infinity, can find a readable introduction in \cite{Tirosh1991}.

\section{Infinity as a limit}
\label{sect:limit}
\subsection{Arithmetic of limits -- A large unspecified number}
In an elementary course of calculus, the first encounter of the students with infinity usually happens as soon as limits of sequences and functions start to be discussed. 

A simple yet useful claim is usually formulated, stating that for any arithmetical operation $\square$ and functions $f,g$, the identity
\begin{equation}
\label{eq:limofarithmeticops}
\lim_{x\to a} (f(x)\square g(x)) = \lim_{x\to a} f(x) \square \lim_{x\to a} g(x)
\end{equation} 
holds true whenever the limits on the right side exist and the operation on the right side is defined. Because a limit of a function can be equal to $\pm \infty$, one now has to extend the usual arithmetical operations $+,-,\cdot, /$ from the set $\R$ of real numbers to the set $\R\suni\lr{\infty, -\infty}$ and, in particular, to specify the cases when an operation remains undefined.

This, of course, can be done definitorically. However no reasonable teacher can expect his students to memorize a list of a few dozen of cases, many of which (as \uvz{$\infty - \infty$ is undefined}, or \uvz{$a/\infty = 0$ for $a\in \R$}) often seem counter-intuitive or nonsensical at the first sight.

A useful metaphor here is to present $\infty$ as a \uvz{very large unspecified number} (and $-\infty$ as a \uvz{very negative unspecified number}). Then all the laws of arithmetics with infinities become fairly intuitive: For example $\infty - \infty$ is indeed undefined, since when one subtracts one very large but unspecified number from another very large unspecified number, the result can be literally anything ($0$ if the two large numbers are identical, a finite number if they are both large but one of them \uvz{a bit} larger than the other, or infinite if one of the large numbers is \uvz{much} larger then the second one). As for the $a/\infty=0$ case, a similar reasoning works --- dividing a reasonably small number by a large number gives a result close to zero, but $\infty$ is \uvz{very large} (larger than any usual number) and therefore $a/\infty$ is closer to zero than any usual number, and, as such, equal to zero.

Naturally, this coarse metaphor should later be replaced by a more correct idea of representing $\infty$ by a function or a sequence that grows beyond all bounds. But such a representation is not possible before a sufficient understanding to the concept of limit is developed.

\section{Infinite cardinalities and ordinalities}
\label{sect:card}

In this section we assume the reader is familiar with basic concepts of set theory, in particular with elementary cardinal and ordinal arithmetics. A good introductory book in the subject is \cite[Chapters 1--3]{JechBook}.

\subsection{Comparing infinite sets -- Shepherd's fair trade}
Quite early in an introductory course of set theory, the students are confronted with the definition of equipotence --- two sets (finite or infinite) $x,y$ are said to have the same size (or to be equipotent) if there is a bijection $f:x\rightarrow y$. This definition, however simple and natural it seems to be to a more experienced mathematician, isn't in any way obvious nor self-explanatory. 
\footnote{In fact, this definition contradicts one of the \uvz{common notions} of Euclid which states that \uvz{the whole is greater than the part}. This discrepancy was a source of a lot of confusion and failed attempts to reason rigorously about infinity throughout the history, and still remains a hard-to-get-over point for students previously exposed only to the classical pre-1900 mathematics.}

One of the common flaws in intuitive reasoning about sizes of infinite sets is to attempt to quantify the size, i.e. to use expressions as \uvz{the number of elements in $x$} and similar ones. This naturally does not make sense until the theory of cardinal numbers is developed. The definition of equipotence is purely qualitative --- it does say when the two sets have the same size without defining what \uvz{the size} is.

In order to explain the difference between \textit{comparing} (qualitatively) and \textit{measuring} (quantitatively) the sizes of two sets, the following example has proven to be effective (see Figure \ref{fig:pasacek}):

\begin{figure}
\includegraphics[width = 0.75\linewidth]{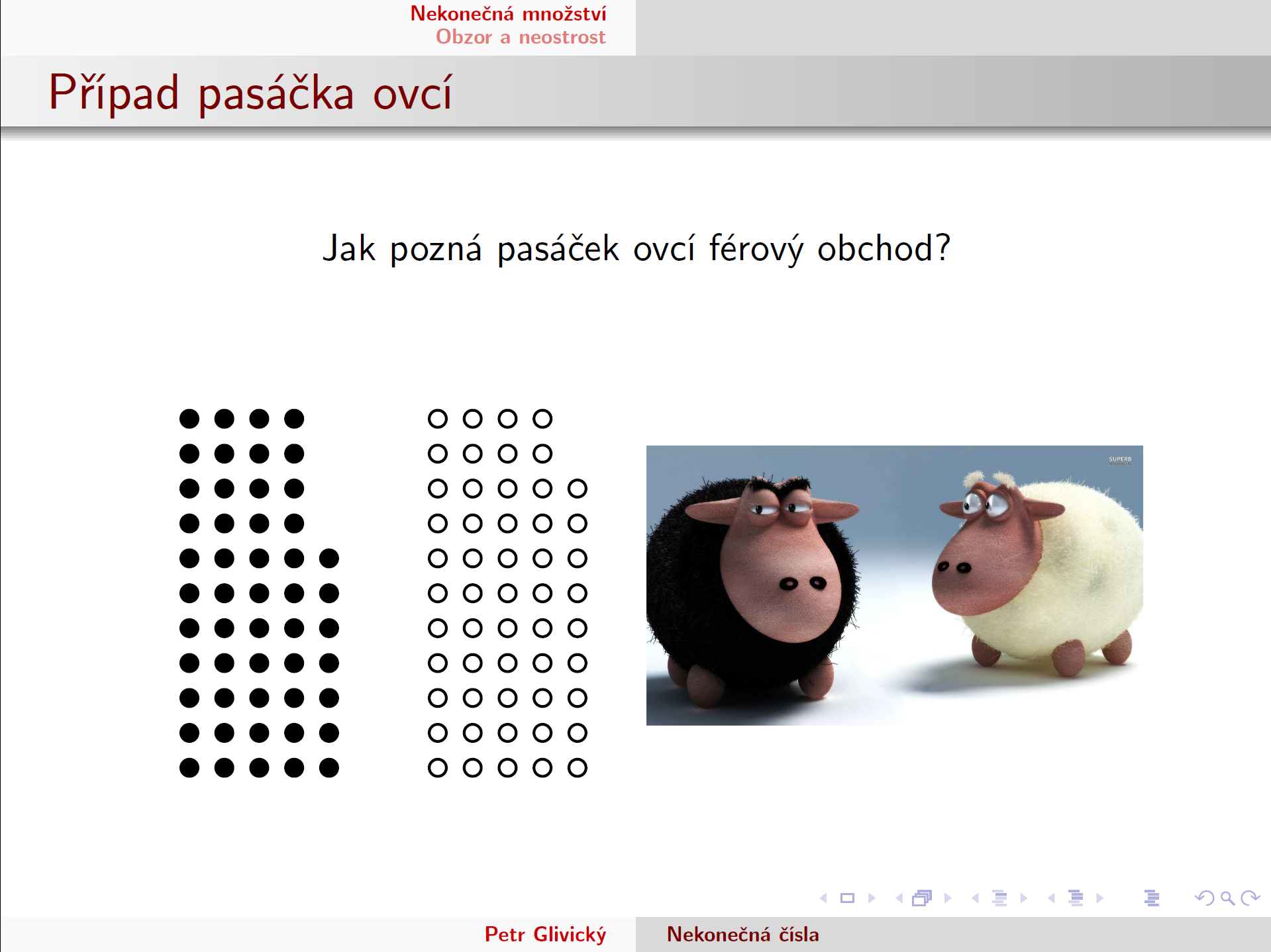}
\caption{A slide on the Shepherd's fair trade example from the presentation \uvz{Infinite numbers, or, do you recognize a monkey?} \cite{MonkeySlides}, in Czech, for high school students. The two flocks of sheep are pictured in such a way that students cannot immediately count their sizes, nevertheless it is evident without counting that there are two more white sheep than the black ones.}
\label{fig:pasacek}
\end{figure}

The students are asked to take a role of an illiterate shepherd owning a flock of white sheep. One day, in the hills, he meets another shepherd with his flock of black sheep. They both like the unusual color of the others' sheep and decide to trade parts of their flocks, a piece for a piece. The black shepherd separates a part of his flock and asks his white counterpart to prepare a flock of the same size for the exchange. If our shephard knew how to count the number of sheep in the offered black flock, he would just count the same number of his whites. But the size of the black flock is much larger than he can count to. Thus he turns out to be in the same situation as we are regarding the infinite sets (flocks) --- we cannot count the number of their elements, because we don't know (yet) such large numbers. The only conceivable way the shepherd can ensure the fair exchange is to exchange the sheep in turns, one for one at a time (thus implicitly constructing a bijection between the exchanged flocks in the process). In the same way we are only able to compare the sizes of the two sets using a bijection but we cannot count their elements.

\subsection{Ordinal and cardinal arithmetic -- Hilbert's hotel}
This is a classical illustration of paradoxical properties of infinite cardinalities, first introduced by David Hilbert in his 1924 lecture \uvz{\"{U}ber das Unendliche} \cite[p.730]{HilbertsHotel}. It is also an excellent example for the basic properties of cardinal and ordinal arithmetics and for the differences between the two. Because this example is rather well-known, we will be little more brief here.

\begin{figure}
\includegraphics[width = 0.75\linewidth]{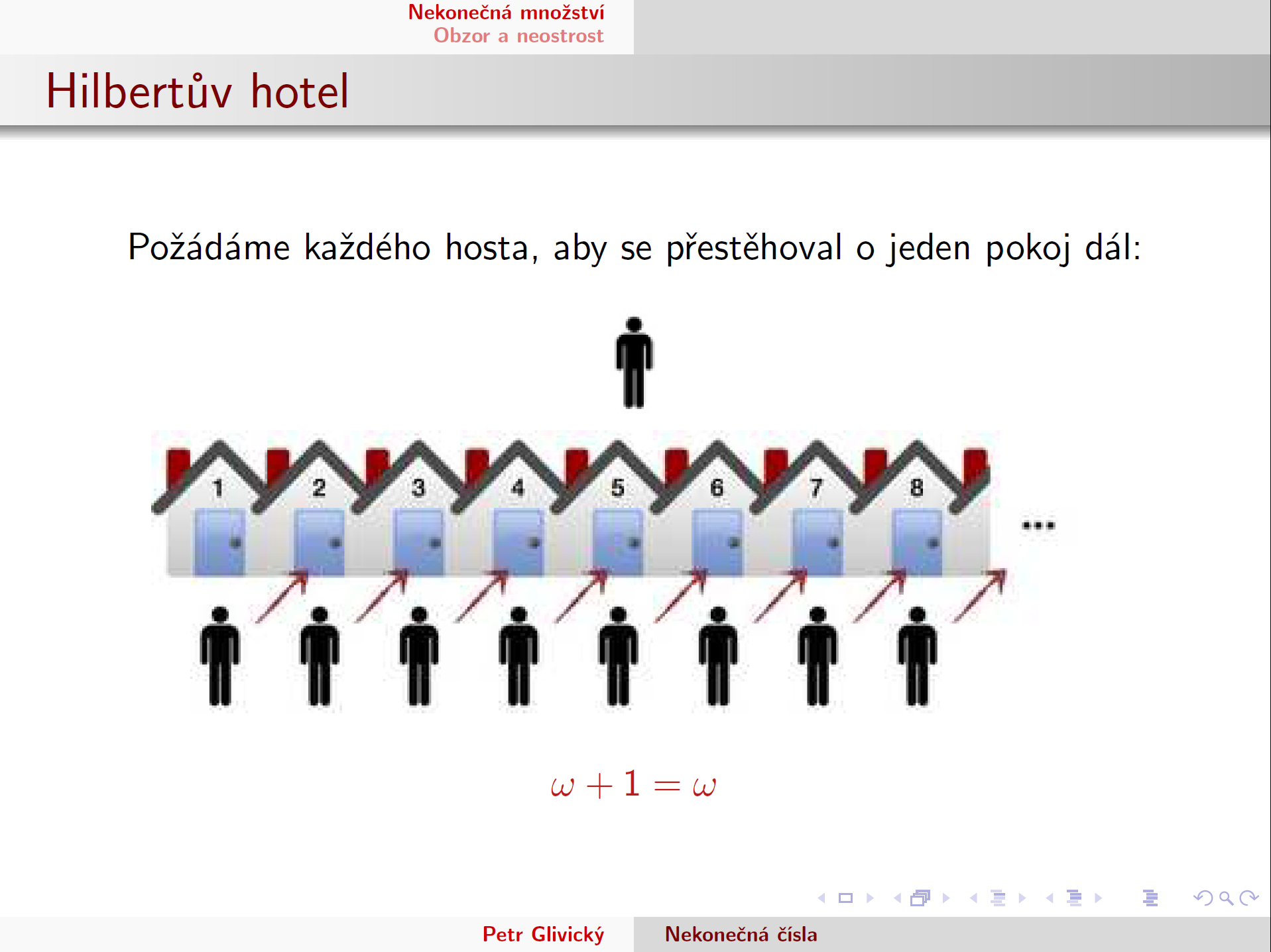}
\caption{A slide on the Hilbert's hotel from the presentation \uvz{Infinite numbers, or, do you recognize a monkey?} \cite{MonkeySlides}, in Czech, for high school students. Here a strategy for squeezing one more guest into a full hotel is represented.}
\label{fig:hilbert}
\end{figure}

The example starts with an idea of an infinite hotel: Let us imagine a hotel with an infinitely (countably) many rooms that are numbered by all natural numbers $0,1,2,\ldots$. Also let us suppose that all the rooms are occupied at the time. The numbering of the rooms is represented by the ordinal number $\omega$ (the order is significant here), while the number of the guests of the hotel can be expressed by $\omega$ as a cardinal number (all the guests are indistinguishable for us). Now several situations (in increasing mathematical profoundness) can be presented: 

First, let us imagine that one new guest comes to the hotel desk and asks for a room (see Figure \ref{fig:hilbert}). Unlike the situation in a fully booked finite hotel, he still can be accommodated --- it is enough to move each guest from the room number $n$ to the room number $n+1$, then the room number $0$ will be empty and prepared for the new guest. This illustrates that in the cardinal arithmetic the equations $1+\omega=\omega+1=\omega$ hold true. 

A more difficult situation to handle is the arrival of an infinite bus (with seats numbered $0,1,2,\ldots$) of new guests into a full hotel. Then, of course, moving each guest from the room number $n$ to the room number $2n$ empties the odd-numbered rooms for the new guests who accommodate them according to the assignment $m\mapsto 2m+1$, where $m$ is a bus-seat number. This illustrates a cardinal equation $\omega + \omega=\omega$.

Yet more complicated situation would be the arrival of the train consisting of countably many train cars, each of which has countably many seats. Accommodating all the passengers into the empty (for the reason of simplicity) hotel is a great illustration of the cardinal equation $\omega\cdot\omega = \omega$. 

Similar stories can be given to illustrate properties of the ordinal arithmetic --- for example a solution to the problem of arrival of one new guest into the full hotel, where we moreover demand that no guest will be forced to move, can be to refurnish a closet at the beginning of the hall to create a new room, to accommodate the new guest there and then renumber all the rooms in increasing order in such a way that the former closet becomes the room number $0$. This illustrates that in ordinal arithmetics it holds $1+\omega=\omega$ (however $\omega+1\neq \omega$ --- if the closet is at the end of the hall, behind all the rooms, then no renumbering that keeps the numbers of the rooms increasing along the hall is possible).

\section{Nonstandard models of arithmetic and analysis}
\label{sect:nonstandard}
We assume here that the reader is familiar with the basic idea of a nonstandard model of Peano arithmetic (see for example \cite{KayeBook} for an elementary introduction to the topic). Conceptual knowledge of the set theoretical framework for nonstandard analysis (see for example \cite{GlivickyNonstandardFramework}) might help to fully appreciate all the nuances of the examples, but the basic ideas can be understood without it.

\subsection{Nonstandard extension of $\N$ -- Blurry horizon}
A fundamental idea that underlies the use of nonstandard methods (such as the nonstandard analysis, or the study of nonstandard models of arithmetic) in mathematics is that of a proper (and often elementary) extension of the structure $\N=\stru{\N}{0,1,+,\cdot,\leq}$ of natural numbers. Such an extension $\starr\N$ of $\N$ is then thought of as the \uvz{true} set of natural numbers, while the original set $\N$ forms a cut on $\starr\N$, i.e. a proper initial segment containing $0$ and $1$ and closed under $+$ and $\cdot$. Such a cut manifests itself as a \uvz{blurry horizon} on the set of natural numbers $\starr\N$ and allows for rigorous definitions of various fuzzy concepts and distinctions as for example infinitesimals, large vs. small, monkey vs. human (see the following sections), and many others.

Of course there is no cut on the set $\N$ (any initial segment containing $0,1$ and closed under $+$ is necessarily equal to $\N$). Considering a \uvz{blurry horizon} on the set of natural numbers thus seems counter-intuitive and is the source of much confusion during the students' first encounters with nonstandard methods. The following metaphor helps to clear things up:

\begin{figure}
\includegraphics[width = 0.75\linewidth]{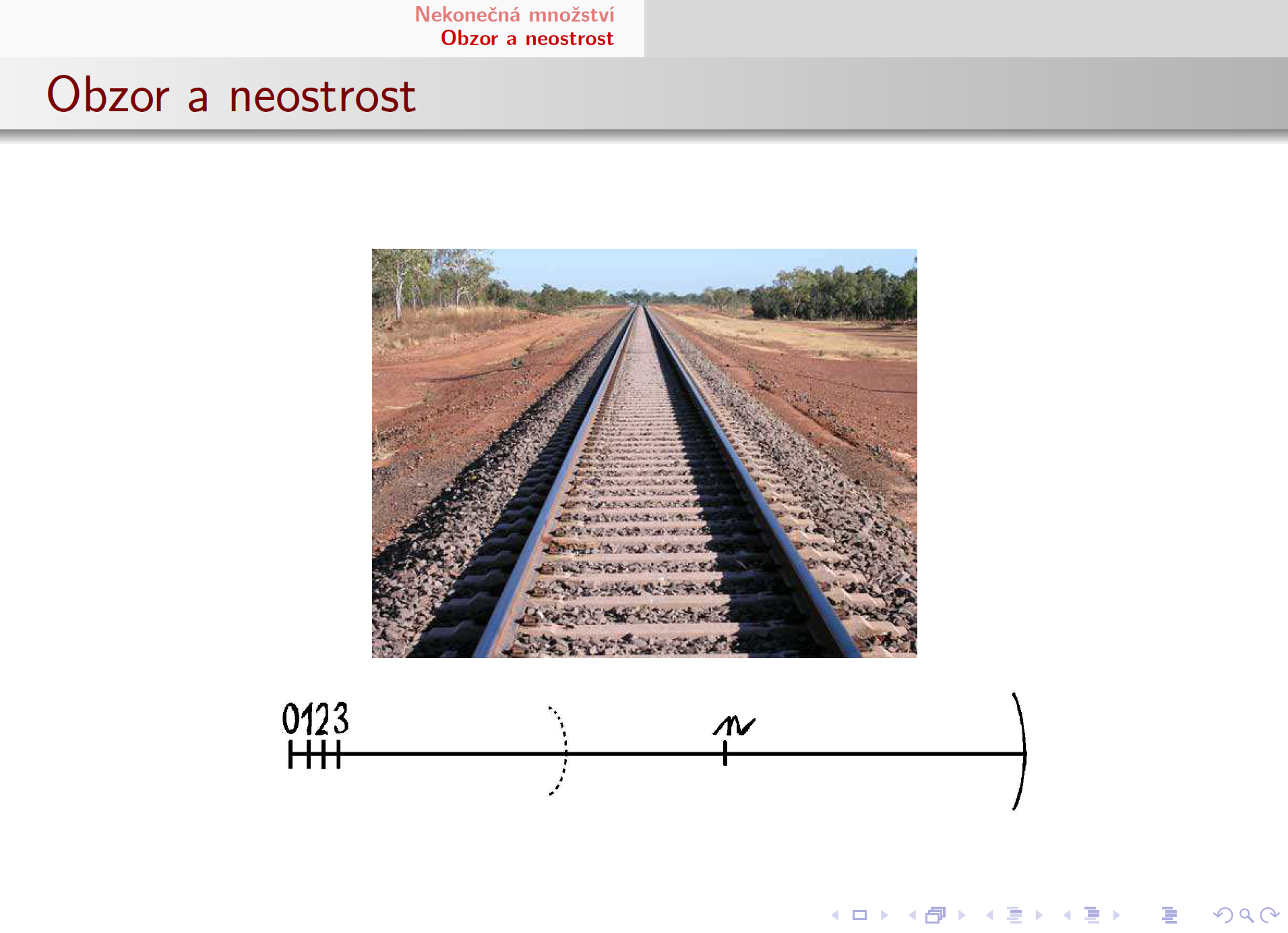}
\caption{A slide on the Blurry horizon from the presentation \uvz{Infinite numbers, or, do you recognize a monkey?} \cite{MonkeySlides}, in Czech, for high school students. The number of railroad ties before the horizon is fuzzy. But clearly there are only finitely many ties between the viewer and any given tie behind the horizon. The same phenomena are manifested in the nonstandard extension $\starr\N$ of $\N$.}
\label{fig:horizon}
\end{figure}

Imagine a railroad stretching from our viewpoint straight to the horizon (see Figure \ref{fig:horizon}) on a perfectly flat ground. In such a situation, classical geometry tells us that we should see the rails as two parallel lines that get seemingly closer to each other as they run further away from us, until they meet in a single point in the infinity. However, in reality, our experience is different. Our eyes are simply not sharp enough to see that far. We stop being able to distinguish the rails long before they get to the point where they meet. Due to imperfection of our sight, the more distant part of the railroad just disappears in a blurry fog. We call this distant limit of our ability to see the \uvz{blurry horizon}.

The whole rail track from our standpoint to the point in the infinity where the rails meet corresponds in our analogy to the nonstandard extension $\starr\N$, while the segment from our standpoint to the blurry horizon represents the cut $\N$ of standard numbers.

Although we cannot see the part of the track beyond the blurry horizon, we know that nothing important changes there --- the railroad looks just the same everywhere. If someone had a better sight than us and was able to actually see the whole track, they would not be able to determine whether they really see more than we do by asking us questions about the track --- their and our experiences would be just indistinguishable (this illustrates that $\starr\N$ is an elementary extension of $\N$).

Let us now focus on the railroad ties (they correspond to individual natural numbers $n\in\starr\N$). When we look at the ties close to us, we can clearly see each one of them individually. Closer to the blurry horizon distinguishing the consecutive ties gets much harder, until, beyond the horizon, we cannot see individual ties at all, however hard we try. Again, the horizon is, in a sense, \uvz{blurry}. Let us try to pinpoint some of the properties of this blurriness.

If we try to count the number of ties we are able to distinguish, or exactly point out the last one that we can see, we never succeed. It almost seems that whenever we catch a sight of one tie, we can also glimpse the next one ($\N$ is indeed closed under taking the successor), but after one or two steps of this process we always blink or let our focus drop a bit and then we need to start over. Despite all this, we know that even for a tie in the blurry beyond-horizon region of our view field, there are only finitely many ties between it and us. 
%Because the rairoad looks the same in front of and beyond the horizon, 
We could, in theory, count exactly how many there are, but we are just not able to keep our focus for a sufficiently long time (numbers of $\starr\N\sdiff\N$ are still natural numbers, but they are \uvz{nonstandard}).
%The horizon is somewhat blurry --- the number of individual ties that we are able to distinguish, seems to change all the time, depending on momentary visibility conditions, our ability to keep focusing, and many other constantly changing conditions. In particular there is never the last tie before the horizon that we can distinguish --- it almost seems that whenever we catch a sight of one, we can also glimpse the next one, but after one or two steps of this process we always blink or let our focus drop a bit and then we need to start over. Despite all this, we know that there are other ties that are beyond the horizon and that we can, in theory,

To make the last idea even more illustrative, imagine that the ties are numbered in the increasing order starting from the point where we are standing. Then the numbers of ties that we can distinguish individually form a \uvz{blurry horizon} on the set of natural numbers. The ties with small numbers lie safely before the blurry horizon, but ties with large numbers are somewhere beyond, in the blur. If we can see the $n$-th tie, we can always (with some difficulties) see the $(n+1)$-th, but it is not humanly possible to iterate this process for a long time, so we can never see ties that are far away.

In the rigorous language of nonstandard methods, the blurry horizon corresponds to the cut $\N$ in the set $\starr\N$. In this analogy, the small natural numbers are all elements of the cut $\N$, but there exist large (called nonstandard) numbers in the \uvz{blurry} part $\starr\N \sdiff \N$. Whenever $n\in\N$, we know that $n+1\in\N$ as well. However by a finite iteration of this process, we can never get outside of the cut $\N$ to the nonstandard part $\starr\N \sdiff \N$.

\subsection{Mathematical induction and nonstandard numbers -- Heap of sand paradox}
The heap of sand paradox was probably first formulated by the ancient Greek philosopher Eubulides of Miletus \cite{Barnes}. In its classical formulation it exploits the vagueness of the word \uvz{heap} in the natural language: 

Let us start with a heap of sand. If we remove just one grain, then what we get, clearly remains a heap. However, since the heap consists only of finitely many grains, by repeatedly removing one grain at a time, we end up with just one grain, which of course is not a heap. So there had to be a moment when by removing a single grain of sand, we changed a heap into a non-heap --- that is paradoxical.

By using the ideas of nonstandard methods, we can not only explain the paradox, but also get more intuition about nonstandard concepts in the process.

\begin{figure}
\includegraphics[width = 0.75\linewidth]{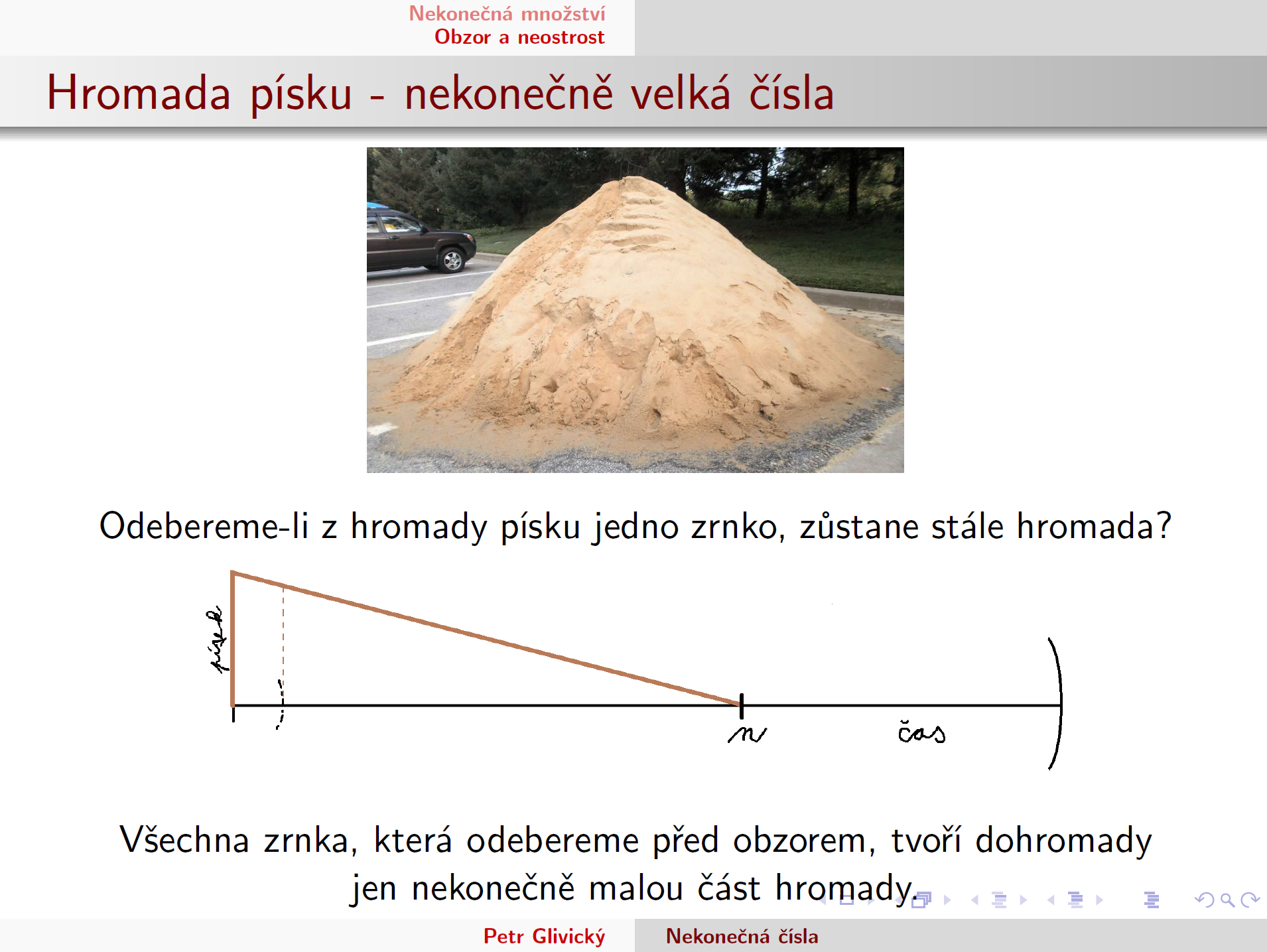}
\caption{A slide on the Heap of sand paradox from the presentation \uvz{Infinite numbers, or, do you recognize a monkey?} \cite{MonkeySlides}, in Czech, for high school students. The fuzzy concept of \uvz{heap} is simulated by using the cut $\N$ in the nonstandard extension $\starr\N$.}
\label{fig:heapofsand}
\end{figure}

We may formulate the paradox mathematically in the nonstandard extension $\starr\N$ of $\N$: The number of grains in the heap of sand at time $t\in\starr\N$ is represented by the function $g(t)=n-t$ for $t\leq n$ and $g(t)=0$ for $t>n$ (see Figure \ref{fig:heapofsand}), where $n\in\starr\N\sdiff\N$ is an initial number of grains in the heap (which is a nonstandard number). The concept of \uvz{being a heap} may be represented as having a nonstandard number of grains (i.e. a number of grains in $\starr\N\sdiff\N$). Then, clearly, if we have a heap at time $t$ (i.e. if $g(t)\in\starr\N\sdiff\N$), then we will have it at time $t+1$ too (as $g(t+1)=g(t)-1\in\starr\N\sdiff\N$). So \uvz{being a heap} cannot be changed by removing one, or finitely many, grains. If we continue in the process of removing grains for a sufficiently long time, we can get to a number of grains in $\N$ (that is to a non-heap), however the time needed for this to happen is nonstandard (i.e. beyond the horizon of our capabilities). Also, it can be seen, that there is no exact moment $t$ where a heap changes to a non-heap --- the transition happens at the blurry horizon $\N$. This is a mathematical manifestation of the fuzziness of the concept of \uvz{heap}.

\subsection{Mathematization of fuzzy concepts -- Charlie the monkey}
\begin{figure}
\includegraphics[width = 0.75\linewidth]{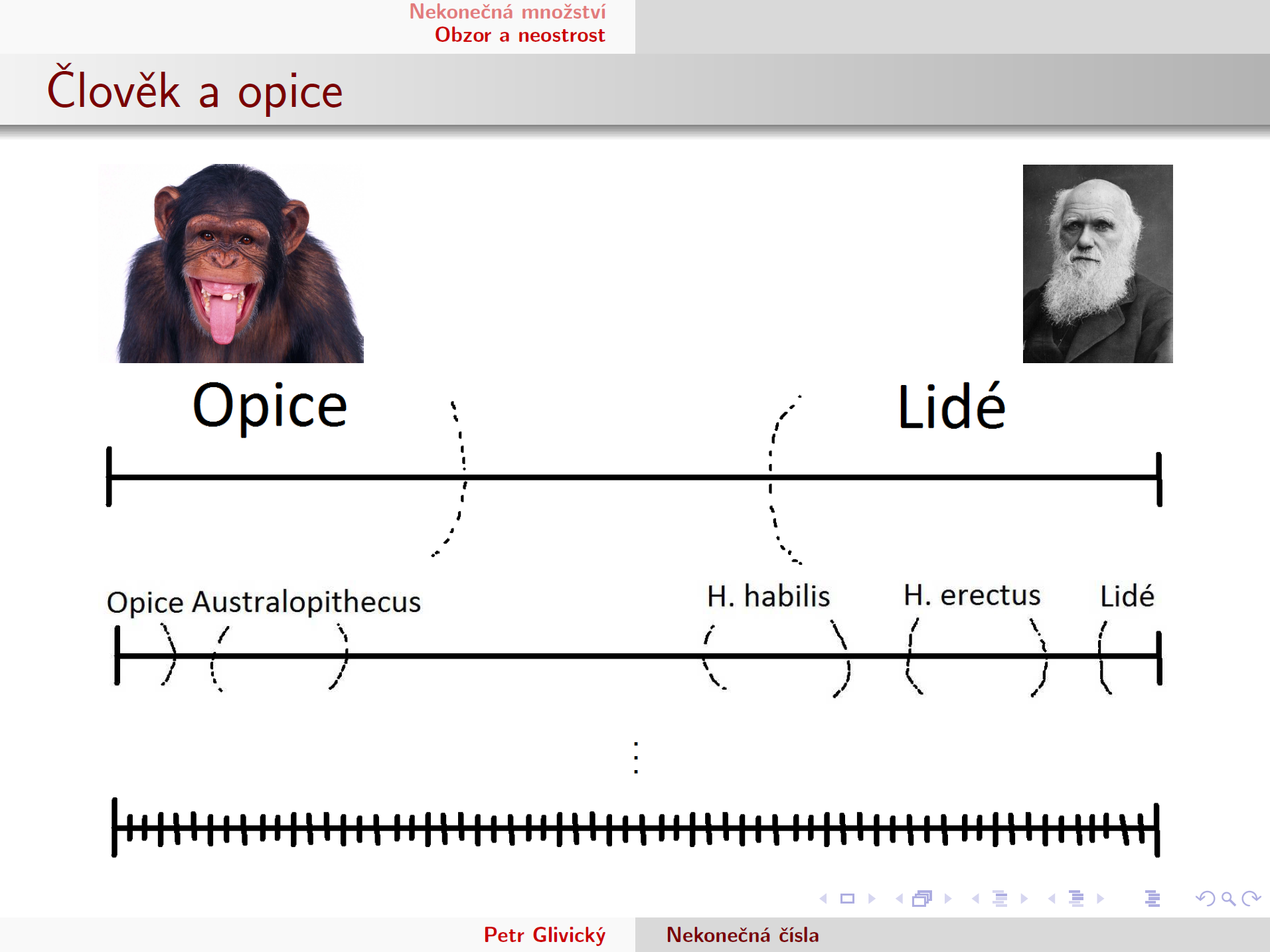}
\caption{A slide on the Charlie the monkey example from the presentation \uvz{Infinite numbers, or, do you recognize a monkey?} \cite{MonkeySlides}, in Czech, for high school students. The fuzzy concepts of \uvz{monkey}, \uvz{human} and many others are simulated by using different cuts on $\starr\N$.}
\label{fig:charlesdarwin}
\end{figure}

The following example is a more complicated version of the previous one. Consider the line of direct male ancestors of Charles Darwin going 350~000 generations back to an ancestor called Charlie the monkey (see Figure \ref{fig:charlesdarwin}), who, in fact, is a monkey. We would agree without hesitation that a father of a human is always a human, while a son of a monkey must certainly be a monkey. But the sequence of descendants from Charlie the monkey to Charles Darwin starts with monkeys and ends with humans. This paradox can be explained by representing the fuzzy concepts of \uvz{monkey} and \uvz{human} in the nonstandard extension $\starr\N$ of $\N$.

For that purpose it is enough to represent the number of generations between Charlie the monkey and Charles Darwin as a nonstandard number $n\in\starr\N\sdiff\N$. Then Charlie is the $0$-th generation, and Charles Darwin the $n$-th one. Being a monkey can then be represented as having a generation number in $\N$, while being a human as having a generation number in $n-\N=\set{n-m}{m\in\N}$. Thus monkeys and humans are two fuzzy concepts at the opposite ends of the ancestral line, leaving a large area of non-humans, non-monkeys in between. As before, we can easily see that there is no last monkey, nor the first human in the line, and that, indeed, the father of a human is a human and the son of a monkey is a monkey.

If we were able to somehow get a picture of each individual in the Charlie-Charles sequence, we wouldn't be able to comprehend them all at once. The number of 350~000 generations is just too big for the human mind (that is why we could represent this number by a nonstandard number $n$). If we tried to see the evolutionary changes in the gradual transition from a monkey to a human, there would be two interesting cases: 
\begin{enumerate}
\item If we took two pictures of individuals only finitely many (i.e. some $m\in\N$) generations away, then the accumulated evolutionary change would be just $m/n$ of the total change, that is an infinitely small part of the total change (as $m\in\N$ and $n\in\starr\N\sdiff\N$). In that case we wouldn't see a difference between the individuals -- the change would be just too small.
\item If, on the other hand, we took pictures of individuals from distant parts of the ancestral line ($m$ generations away, where $m$ is approximately equal to $n/k$ for some $k\in\N$), then the accumulated change $m/n$ would be approximately $1/k$ of the total change, which for small $k$ is quite noticeable. In that case we would be able to identify the two individuals as specimens belonging to two different species of ancestors of modern humans.
\end{enumerate}

\bibliography{biblio}{}
\bibliographystyle{amsalpha}

\end{document}